\documentclass[twoside,11pt]{article}                                             %
\usepackage{graphicx}
\newtheorem{thm}{Theorem}
\newtheorem{corr}{Corollary}
\newtheorem{lemma}{Lemma}

\newtheorem{prop}{Proposition}
\newcommand{\proc}{\textbf}
\newcommand{\ep}{}

\begin{document}
\font \cmbxslten=cmbxsl10 
\def\La {\mbox{\cmbxslten \char'003}}
\font \cm=cmbxsl10 
\font \cma=cmbxsl10 scaled\magstep1
\font \cmbten=cmb10 
\def\Ga {\mbox{\cmbxslten \char'000}}
\font \cba=cmb10 
\font \cbb=cmb10 scaled\magstep1 \font \cbc=cmb10 scaled\magstep2
\font \cbd=cmb10 scaled\magstep3
\pagestyle{myheadings}                
\markboth{\hfill Alexandr Prishlyak} {Topological classification
of Morse--Smale flows on 3-manifolds}
\thispagestyle {empty}                                                   %
\setcounter {page} {1}                









\title{Topological classification of Morse--Smale
flows on 3-manifolds}

\author{ALEXANDR PRISHLYAK }


\maketitle

\begin{abstract}
  We construct  a topological invariant  for a
  Morse--Smale flow on a  $3$-manifold and prove
  that the flows are  topologically equivalent iff
  their  invariants are same.
\end{abstract}







\section{Introduction}

In this paper we study the topological properties of Morse-Smale
flows on closed oriented 3-manifolds. We give the topological
classification of these flows by constructing complete topological
invariant (distinguished graph) and proving the criteria of
topological equivalence.

Let $M$ and $N$ are closed smooth manifold.

Two flows $X$ and $Y$ on $M$ are topologically equivalent, if
there exists a homeomorphism $h: M \to M$  that carries
trajectories of $X$ to those of $Y$ and preserves their direction.

 In \cite{P73}, Peixoto
gave  a topological classification of Morse--Smale flows on closed
orientable 2-manifolds.  There are many papers where this
classification was improved and  new approaches were proposed. For
references, see, for example \cite{AG86,F75,F83,OSh98,P03,W90}.

Vlasenko \cite{V} and, independently, Bonatti and Langevin
\cite{BL98} gave a classification of Morse-Smale diffeomorphisms
on closed 2-manifolds.

A topological classification of Morse--Smale dynamical systems on
3-manifolds are known only for some classes of systems
\cite{BB01,F75,P98,P02,R87,U90}.










For Morse--Smale flows without closed orbits, Smale constructed
handle decompositions. For flows with  closed orbits, Asimov
constructed a round  handle decomposition. For our classification
of Morse--Smale flows, we use Asimov's idea that it is possible to
reconstruct the flow by using the decomposition on prime and round
handles. However, it is necessary to note that there are
infinitely many such decompositions for the same flow. Therefore,
it is necessary to chose among them  a canonical one. In the
general case, a Morse--Smale flow  has infinitely many
trajectories in  the intersection of 2-dimension stable and
unstable manifolds. This difficulties are similar to difficulties
that arise in the classification of Morse--Smale diffeomorphisms
on surfaces. Besides,  for dimension $3$,  it is necessary to
introduce an additional invariant  for an extension of
homeomorphisms from the boundary of the round handles to their
interior. We call this invariant a $\tau$-invariant.

\section{Basic definitions}

Let $M$ be a closed 3-manifold and $X $  a flow on it. $X $ is  a
Morse--Smale flow, if the following conditions hold:

1) $X $ has a finite number of critical elements (fixed points and
closed orbits) and all of  them are non-degenerate (hyperbolic);

2) stable and unstable manifolds of critical elements have
transversal intersections;

3) the non-wandering set of $X$ is a union of fixed points and
closed orbits.


Two flows $X$ and $Y$ on $M$ are topologically equivalent, if
there exists a homeomorphism $h: M \to M$  that carries
trajectories of $X$ to those of $Y$ and preserves their direction.

Denote by $v ( x )$ and $u ( x )$ stable and unstable manifolds of
a critical element $x$. If $x $ and $y $ are  critical elements
such that $\dim v(x) = \dim u (y) = 2$ and $v(x) \cap u(y) \ne
\emptyset $, then [7] there is a sequence of saddle closed orbits
$ \beta _{1}, ..., \beta_n$ satisfying $$v (x) \cap u(\beta _{1} )
\ne \emptyset, v(\beta _{k} ) \cap u(b _{k +1}) \ne \emptyset
,v(\beta_n)\cap u( y ) \ne \emptyset .$$

Denote by $m $  the maximal number of elements $ \beta _{k} $  in
all such sequences and  set {beh}\,$(y \vert x) = m + 1$. Thus
{beh}\,$( y \vert x) = 0$, if $v(x) \cap u ( y) = \emptyset $. We
say that a flow has beh equal $s$, if $s $ is a maximal value of
beh\,$(y \vert x)$.

The number beh\,$(y) = \max _{x}$ beh\,$( y , x )$ is called  the
height of the closed orbit of  index 1. Thus, a closed  orbit  has
height 0, if its stable manifolds do not intersect stable
manifolds of closed orbits of  index 1.  A closed orbit of index 1
has height 1, if its stable manifolds intersect unstable manifolds
of closed orbits of height 1 and do not intersect unstable
manifolds of other closed orbits of
 index 1, etc.

Bellow we assume that $X $ and $X ' $ are Morse--Smale flows on a
closed oriented 3-manifold $M$.

Let $a _{{ 1}}..., a _{k} $ and $a _{{ 1}} '..., a _{k} ' $ be
index 0 fixed points of of $X$ and $X ' $, respectively, and let
$b _{{
    1}}..., b _{n} $ and $b _{{ 1}} '..., b _{n} ' $ be index 1 fixed
points of $X $ and $X ' $, correspondingly. By $ \gamma _{j} ^{i}
, \ j = 1, \ldots, m(i)$, we denote closed orbits the height of
which equals $i$.

Let $K$  be the union of  the following stable manifolds:

1) fixed points of  indices $0$ and $1$,

2) closed orbits of  index $0$.

Fixed points of  index $1$ are called lower points, and their
unstable manifolds are called  lower manifolds.

Consider a regular neighbourhood $U (K) $ of $K $. Let $\Phi =
\partial U (K)$  be the boundary of this neighbourhood of the flow
$ X $ and  let $\Phi ' $ denote the boundary for the flow $X'$.

Fixed points of  index 2 and closed orbits of  index 1 are called
upper critical elements; their stable manifolds are called upper
manifolds. The curves that are intersections of $\Phi $ and lower
manifolds are called lower curves, and the intersection of $\Phi $
and upper manifolds are called upper curves.

Note that all lower curves are closed (homeomorphic to circles),
and the upper curves can be both closed and not closed. Moreover,
infinitely many open curves can correspond to an upper critical
element.

The surface $\Phi $, together with the lower and upper curves,  is
called a diagram of the flow. Two diagrams are called
homeomorphic, if there exists a homeomorphism of $\Phi $ that maps
the lower and upper curves to lower and upper curves,
respectively.

\begin{prop}  If flows are topologically
equivalent, then their diagrams are homeomorphic.
\end{prop}

\proc{Proof.} Let $ \psi $: $M \to M ' $  be a topological
equivalence between $X $  and $X ' $. Then the required
homeomorphism of the diagrams is  given by the formula
$$\Phi \mathrel
x \to \varphi  ( \psi  (x ))  \cap \Phi '
,$$ where $ \varphi  (y )$ is  a trajectory of $X ' $ that
 contains $y. $
\ep \medbreak

If the flow doesn't  have closed  orbits, the existence of a
homeomorphism  of the diagrams is also a sufficient condition  for
topological equivalence of  the flows \cite{P02}.

This is not true in general and, in addition,  it is  difficult to
prove  that there exists a homeomorphism  between diagrams with an
infinite number of noncompact curves.

In what follows, our aim  would be to  apply the diagram to the
surface $\Phi $ with the embedded graph $G $ and  determine
additional information which, together with  the pair ($\Phi $, $G
$), would  give a complete invariant of  the flow.  It is
convenient to express this information  in terms of the
decomposition  in simple and round handles.

A handle decomposition of $M$ is a sequence $D^3=M_0 \subset M_1
\subset M_2 \subset ... \subset M_k = M$, where $M_i$ is obtained
from $M_{i-1}$ by adding a handle or a round handle.

We construct a handle decomposition such that the cores and the
cocores lie on stable and unstable manifolds of the critical
elements, respectively.  We order the critical elements in such a
way that, for all $i $, the following holds: critical points of
index $i $ are less than closed orbits of index $i $ which are
still less than critical points of index $i $ + 1; for closed
orbits of index 1, $\gamma _{j} ^{i} < \gamma _{j} ^{i{ +}{ 1}}$.
Denote the critical elements by $x _{{ 0}}, x _{{ 1}} $, \ldots,
$x _{n} $ ($x _{{ 0}} < x _{{ 1}} < \ldots < x _{n} $).  Then we
construct a decomposition in handles in the following way: the
handle  $v(x_0)$ is a closed tubular neighbourhood $U (x_0)$ of
the least critical element $x _0$ , $H _{i} = $cl$(U (v(x _{i} ))
\backslash \cup _{k{ <} i} \quad U _{k} ), \ i = 1, \ldots, n $.

Note that the surface $\Phi $ can be considered as the boundary of
the union of handles that correspond to critical elements of index
0 and critical points of index 1.  Since stable and unstable
manifolds intersect transversally, it is possible to choose the
neighbourhood $U $ ($ v $ ($ x _{i} $)) to be so small that all
the intersections of the core and the cocore with the bases of
$2$-handles and round $1$-handles, as well as the side legs of the
$1$-handles and the round 1-handles, are parallel.  This means
that it is possible to introduce a structure of the direct
products, $S ^{{ 1}} \times [-1,1]$, such that these intersections
have the form $\{s _{i} \} \times [-1,1]$.

If the flow  does not have closed orbits, then all  decompositions
in handles, which have parallel  intersections, are isomorphic and
can be given, up to an isomorphism, by the diagram. In case of a
flow with closed orbits of height not smaller than 2, a decrease
of neighbourhoods of their stable manifolds  increases the number
of intersections of the cores and cocores of the handles of
smaller  orders. The idea of our construction is to  choose
handles at each step  to be the largest neighbourhoods for which
the property of the intersections is preserved. If the constructed
handle decompositions are isomorphic, then the complexes  that
consist of the boundaries of the handles are isomorphic and, also,
the isomorphism maps meridians of the round handles into
meridians.

As  in Proposition~1,  different choices of simple 0- and
1-handles and round 0-handles  lead to isomorphic handles
decompositions. We consider in more  details the case of round
1-handles.

\section{ The structure of the flow  in a neighbourhood
  of a closed 1-orbits}

Closed orbits of  index 1 can have trivial or twisted
neighbourhoods (both of them are solid  tori). In case of a
trivial neighbourhood, it intersects stable and unstable manifolds
in cylinders $S ^{{ 1}} \times I $ and  the torus  in a pair of
circles. In case of twisted neighbourhoods, their intersections
are Mobius bands and circles.

A trivial neighbourhood can be representaed as the product $[1,3]
\times S ^1 \times [-1,1]$ with coordinates $\{ \rho , \alpha , z
\}$ (an analog of the  cylindrical coordinates) \cite{P98}. Thus,
if a trajectory intersects a torus in a point with $ \rho =2 \pm
1, \alpha = \alpha _0, z=z _0 \ne 0$, then its second intersection
point with the torus  has the coordinates
\begin{equation}
\rho = 2 \pm \vert  z _0 \vert, \ \alpha  = \alpha _0 + \ln \vert
z _0 \vert , \ z =  \mbox{sign } z_0.
\end{equation}

A stable manifold is  determined by the equation $z = 0$,  an
unstable manifold by $ \rho =2 $.

We consider the torus $T $ as the union $T=T^- \cup T^{+} $, where
$$ T^{{ -}}  = \{ \rho \in  \{1; 3 \}, -1
 1 \},$$
$$ T ^{+}  = \{1
3, z \in \{-1; 1 \} \}.$$

Thus, $T ^{{ -}} \subset \Phi $ is a set of the points incoming in
torus  and $T ^{{ +}} \subset F$ is a set of outgoing points.

Formula (1)  defines a map $g $  that maps points of $T ^{{ -}} $
with coordinates $ \rho =1, z > 0$ into points of $T ^{{ +}} $.
Similar formulas
 hold if $ \rho =3$ and $z < 0$.

The intersection of the stable manifold of a closed 1-orbit and
the torus consists of two circles, $ \rho =1, \ z = 0$ and $ \rho
=3, \ z = 0,$ which we denote by $v ^{{ 1}} $, $ v ^{{ 2}} $. The
unstable manifold of  a closed 1-orbit intersects the torus  in
two circles, $ \rho =2, \ z =- 1$ and $ \rho =2, \ z = 1,$ denoted
by $u ^{{
    1}} $, $u ^{{ 2}} $




Since stable and unstable manifolds have transversal
intersections, the curves $u _ {i} $ intersect the circles $v ^{{
1}} $, $ v ^{{ 2}} $ also transversally.  Thus, in a neighbourhood
of an intersection point, $u _{i} \cap v ^{{ 1}} $, the curve $u
_{i} $ is  given by the equation $ \rho $ =1, z = $ t, \ \alpha =
t k (t),$ where $k=k (t) $ is a smooth function. The image of this
curve  with respect to the map $g $ is
 the
curve 

$$ x = (2- t ) \cos (t k (t) + \ln t ),$$ $$y = (2 -t ) \sin (t k
(t) + \ln t ),$$ $$z = 1.$$



 As $t \to 0 +$, this curve  winds
  around $u ^{1} = \{ \rho =2, z = 1 \}$, and the tangent to the
curve  approaches the tangent to the circle. As $t \to 0-$,  the
same is true if
   $z = 1$ is replaced by
$z = - 1$ and $t $ by $ \vert \quad t \vert $.

Similarly, for arcs in the intersections $u _{i} \cap v ^{ 2} $
with $ \rho $ =3, the image of the arc of $u _{i} $  winds around
the same circles $u ^{{ 1}} = $ \{$ \rho $ =2, $z = $ 1 \} and $u
^{{ 2}} = $ \{$ \rho $ =2, $z =- $ 1 \}.

For the intersections $v _{i} \cap u ^{{ 1}} $ and $v _{i} \cap u
^ {{ 2}} $ from the upper and lower  bases, the cluster sets of
the pre-image $g ^{-{ 1}} $ ($v _{i} $)
 will contain the circles $v ^ {{ 1}} $ and $v ^{{
    2}} $.

By reducing the neighbourhood of the closed orbit, it is possible
to achieve that  the following two properties hold:

1) Each arc  in the intersection $u _{i} \cap T $
 intersects $v ^ {{ 1}} $ or $v ^{{ 2}} $
transversally in one point.  Each arc  in the intersections $v
_{i} \cap T $  intersects $u ^{{ 1}}$ or $u ^{{ 2}} $
transversally in one point.

2) If an orientation is fixed on each arc  in the intersections $u
_{i} \cap T $ and $v _{j} \cap T $ and  on the torus $T $, then
all the intersections $g $ ($ u _{i} $) $ \cap v _{j} $ have  the
same sign of intersection. Similarly, all the intersections $g
^{-{ 1}} (v _{j} ) \cap u _{i} $ also have the same sign.

Here the sign of the intersection is positive, if the orientations
of the  intersecting curves determine an orientation of the torus
such that it coincides with the given one; the sign is negative in
the opposite case.

A neighbourhood of the closed orbit satisfying these properties is
called standard.

For  a twisted neighbourhood, the concept of a standard
neighbourhood is introduced  in a similar way. It can also be
obtained from a standard trivial neighbourhood by slitting it into
the disks $ \alpha $=const and subsequent pastings a pair of the
obtained disks together via the central symmetry map 


Let us consider  these neighbourhoods as round 1-handles.   The
orientation of the closed orbit  introduces a parallel orientation
in the corner $T^- \cap T^+$.

\section{The scheme of the  flow }

Let us consider the diagram of a flow. On the lower and upper
circles and the arcs that correspond to  closed orbits, we set an
orientation to be parallel to the direction in which the
corresponding closed orbit moves. We pair these oriented circles
and sets of arcs together  in such a way  that the circles or the
sets of the arcs from one pair correspond to  a  single closed
orbits.  Since the set of arcs can be considered as  cuts of the
circle, there is a natural cyclic order on them.

We consider the upper circle $u _{i} ^{0}$ which corresponds to
the closed orbit $ \gamma _{i} ^{ 0} $ of  height 0. From the
previous paragraph it follows that there is  a sufficiently small
neighbourhood $U $ of the circles $u _{i} ^{{ 0}} $ on $\Phi $
satisfying  the  following property: it can be cut  along a curve
that is transversal  to all upper circles and arcs so that the
obtained set is homeomorphic to $[0,1] \times [0,1]$, and the cut
of the upper arcs in $U $  has the form $\{t \} \times [0,1]$.  A
neighbourhood with such properties is called trivial. If $u _{i}
^{{ 0}} $ intersects the lower circles, then we take the arcs of
these circles  to  be a  transversal curve.  We demand that the
boundary of $U $ consist of  arcs of the diagram. A trivial
neighbourhood $U $ is called  maximal (MTN), if there  are no
greater neighbourhoods  that possess the  above properties 
Fig. 4).



For circles of height more than 0, a maximal  trivial
neighbourhood is defined similarly. It is easy to see that  a
maximal neighbourhood always exists, except for  the case when
there is an annulus $S^{1} \times [0,1]$ such that its boundary
consists of upper circles and all the arcs are wound around $S^1
\times \{0 \}$ with one end and around $S ^1 \times \{1 \}$ with
another end (see Fig. 5). In this case, we take the maximal
neighborhood to be any trivial not intersecting neighborhood, the
complement to which has in $S ^{{ 1}} \times [0,1]$ the minimal
number of points of intersection of the lower circles, the upper
circles, and the arcs.


Consider now the set of arcs $\{u \}$ whose diagrams correspond to
a closed orbit of height $1$. Since its stable manifold intersects
unstable manifolds of dimension $2$ in a finite number of
trajectories, this set consists of a finite number of arcs.  Let
$U $ be a neighbourhood of the set of arcs on $\Phi $ in the
complement of a maximal trivial neighborhood of circles of height
$0$.  Such a neighbourhood is the union of products $[-1,1] \times
[0,1]$. It is called trivial, if all its intersections with arcs
of height not less $1$ are of the form $\{t \} \times [0,1]$.  The
boundary is also considered as consisting of arcs. Besides, we
require that all the arcs $\{ \pm 1 \} \times [0,1]$ correspond to
a single closed orbit, and  the the relation of cyclic order on
them argees with the relation on $ \{u \}$.  A trivial
neighbourhood is called maximal if there is no a greater maximal
neighbourhood possessing these properties.

Similarly to the case of circles, there may not exist a maximal
neighbourhood. In such a case, the maximal neighbourhood is taken
to be an  arbitrary neighborhood whose complement does not contain
points of intersection of arcs and circles.

Consider a collection of arcs $\{w \}$ that correspond to a closed
orbit of height $2$. If the arcs of height $1$ are  considered as
a circle cut in a finite number of points, then $\{w \}$ can be
considered as a finite number of arcs winded around these circles.
Cutting the circle with a transversal curve leads to cutting the
winding arc into an infinite number of arcs. Hence, $\{w \}$
consists of an infinite number of arcs. Nevertheless, there is
only a finite number of them that do not lie in the maximal
neighbourhoods constructed above.  For these arcs, a construction
of a maximum trivial neighbourhood is same as for arcs of height
$1$.

Again, there is only a finite number of arcs of height $3$, not
lying in maximal neighbourhoods, and arcs of heights $0$--$2$. For
them, we construct maximal neighbourhoods, and we continue this
process until maximal neighbourhoods are constructed for all sets
of the arcs.

The constructed maximal neighbourhoods are used for constructing
round handles. These neighbourhoods will be intersections of the
corresponding round handles and $\Phi $. We construct the round
handles starting with height $0$.  Round handles of height $0$ are
neighbourhoods of stable manifolds in the way that their bases are
the maximal neighbourhoods. The arbitrariness in the choice of the
side walls is not important, since the first return map defines a
homeomorphism between them (and our constructions are, actually,
carried out up to a homeomorphism).  Since the lower bases are
trivial neighbourhoods, it follows from Section 2 that the
intersections of stable manifolds of height $0$ and side walls
make parallel curves. In side walls of $1$-handles, consider
regular neighbourhoods of parallel curves of height $1$, which
extend maximal trivial neighbourhoods of height $1$.  This means
that these neighbourhoods contain parallel curves of height more
than $1$ and that their ends lie in MTN of height $1$. The union
of these neighbourhoods and   the MTN of height $1$ are bases of
round handles of height $1$. Similarly to the case of handles of
height $0$, various choices of side walls lead to homeomorphic
constructions. The stable manifold of height more $1$, again,
intersects these side walls in parallel curves.  In the
constructed side walls of heights $0$ and $1$, we choose
neighbourhoods of parallel curves of height $2$ that  would extend
the MTN of height $2$ to the bases of round handles of height $2$;
round handles of all other heights are constructed  similarly.

We thus have constructed some two-dimensional compact stratified
set $S $, which is the union of the surface $\Phi $ and the
boundaries of round handles of index $1$.  In this construction,
the union of $0$- and $1$-dimensional strata consists of the
following closed curves: 1) b-spheres of $1$-handles, 2) a- and
b-spheres of round $1$-handles, 3) a-sphere of $2$-handles, 4)
boundaries of the bases of round $1$-handles.  On strata of the
second type, the orientation is defined in correspondence with the
motion on the corresponding closed orbit of index $1$.  It also
follows from  the construction that the diagram of the flow can be
recovered uniquely (up to a homeomorphism) from the stratified
set.  Thus, existence of an isomorphism of the stratified sets is
a necessary condition for topological equivalence of the flows.

On each torus, which is the boundary of a round handle, we single
out a closed curve in $S $, which bounds a disk that intersects
the corresponding closed orbit transversely in a single point.

\proc{Definition.} \textit{A scheme} of the flow is a stratified
set $S $, together with the following: 1) the set of
$2$-dimensional subsets that are boundaries of round handles and
their indices, 2) a partition of $1$-strata into 4 types, 3)
orientations of cycles of the 2nd type and partitions of some of
them into pairs.

Two schemes are called isomorphic if there is a homeomorphism of
the stratified sets, which preserves types of the $1$-stratas and
the partition of the cycles into pairs and orientations (on those
cycles where the orientations are given).  \medbreak

\begin{thm}  There is a topological equivalence of
two Morse--Smale flows $X$ and $X'$ preserving the order of
critical elements if and only if there is an isomorphism between
their schemes mapping chosen curves of the first scheme  into
curves that are homotopic to the chosen curves of the second
scheme and preserving the chosen orientations on these curves.
\end{thm}

\proc{Proof.} \textit{Necessity.} Let flows $X $ and $X '$ be
topologically equivalent. Then there is a homeomorphism $h$ that
maps trajectories of $X $ into trajectories of $X ' $. The image
of the scheme of $X $ is a scheme of $X' $. Thus, it is necessary
to show that two different schemes of the same flow are
isomorphic. Indeed, the first return map defines a homeomorphism
between the diagrams and, also, a homeomorphism between the
stratified sets. Thus, since another choice of  chosen curves
yields homotopic curves by the construction, we  see that the
first return map is an isomorphism between the schemes.

\textit{Sufficiency.} Let an isomorphism between the schemes of
Morse--Smale flows be given. Let   us show how to construct
another isomorphism between the schemes such that the
homeomorphism $f : S \to S' $ can be extended to a homeomorphism
of the manifolds. By using the construction and the fact that the
bases of round $1$-handles are trivial neighbourhoods, it is
possible to restore the diagram from the scheme, and the scheme
isomorphism induces a diagram isomorphism $\Phi \to \Phi $.

\textit{The first step.} By given isomorphism of the scheme we fix
homeomorphisms of angles of the round 1-handle. Consider round
handles $H _{i} ^{ k} $ of height $k$ for closed orbit $ \gamma
_{i} ^{ k} $. Stable and unstable 2- dimensional manifolds
intersect the boundary of the handle in curves of two types: 1)
transversal intersected stable or unstable manifold in one point,
2) noncompact. Consider  the first intersection map $g $ from the
bases to the side walls. It map curves of one type to the curves
of other type (see Fig. 6). The isomorphism of the scheme give one
to one correspondence between the intersections of the curves. We
extend it to homeomorphisms of noncompact curves of the bases of
$H _{i} ^{ 0} $ coordinated with homeomorphisms of the angles.
Using the first intersection map $g $ we construct the
homeomorphisms of transversal curves  $H _{i} ^{ 0} $, then we
extend it to homeomorphisms of nonocompact curves of the bases of
$H _{i} ^{ 1} $ and so on. By analogy, we construct homeomorphisms
of noncompact curves of side walls and transversal curves of bases
beginning from the handles of biggest height. We extend the
homeomomorphisms of curves to regions bounded by them on bases of
$H _{i} ^{ 0}$. Using $g$ construct homeonorphisms of region of
side walls of $H _{i} ^{ 0}$. The construct homeomorphisms of
regions of f $H _{i} ^{ 1} $ and so on. Thus we construct the
homeomorphisms of the boundaries of round 1-handles.


Since the chosen  cycles are homotopic we can extend the
constructed homeomorhism to the disk having one transversal
intersection with the closed orbit and chosen cycles as it
boundary.

Now it is possible to extend the homeomorphism of tori to
homeomorphisms of solid tori.

 The latter is equivalent to the
condition that each circle of the torus $T _{i} $, which bounds a
disk in the solid torus, is mapped into a circle that has the same
property. Then we can extend the homeomorphism of the circle to a
homeomorphism of the disk so that it will be compatible with the
first return map from the torus to the disk.  Furthermore, the
complement of this disk in the solid torus is a three-dimensional
disk and, using compatibility of the homeomorphism of the torus
and the map $g $, we can extend the homeomorphism from the
boundary of the $3$-disk inside it along the trajectories.  To do
this, fix a Riemannian metric on the solid torus. We consider
parts of the trajectories that lie in the $3$-disk and have their
ends in its boundary. Then, by the construction of the
homeomorphism on boundary, it maps each pair of ends of one
trajectory into a similar pair. We define a homeomorphism between
these parts of the trajectories in such a way that all arcs of
these trajectories have proportional lengths.  So the constructed
homeomorphism of parts of the trajectories defines a homeomorphism
of the solid tori.

\textit{The second step.} Let us extend the homeomorphism from
$\Phi $ to $U (K) $.

For each lower closed orbit $ \alpha _{i} $ of $X $, we choose a
neighbourhood as in Section 2 and a sufficiently small
neighbourhood of $u _{i} $ satisfying $T _{i} ^{+} \subset \Phi $.
For the flow $X ' $, these neighbourhoods are chosen so that ${T
_{i}^{+}}' = f ( T _{i} ^{+} ).$

Then there exist natural homeomorphisms of such neighbourhoods so
that they coincide with $f $ on $T _{i} $. Since $ \alpha _{i} $
is a lower closed orbit, we can assume that $T _{i} ^{+} $ does
not intersect $u _{j} $ except for $u _{i} ^{{ 1}} $ and $u _{i}
^{{ 2}} $.

For each fixed point of index $1$ there is a neighbourhood which
is homeomorphic to the cylinder $\{ \rho
1 \}$. Assume that the side walls of such cylinders
are mapped into sufficiently small neighbourhoods of the curves $u
_{i} $ on $\Phi $. For the flow $X ' $, let us construct cylinders
so that their side walls coincide with the images of the side
walls under the mapping $f $. There are natural homeomorphisms
between the cylinders, which coincide with $f $ on the side
surfaces.

If we remove, from $U (K) $, the constructed cylinders, which are
neighbourhoods of fixed points of index $1$, and solid toruses,
which are neighbourhoods of the lower closed orbits, we obtain a
union of $3$-disks each of which contains one source. Having a
homeomorphisms of the boundaries of these disks, we extend it
inside along the trajectories.  We thus obtain the required
homeomorphism $U (K) $.

Extend the homeomorphism to the rest of the manifold $K$ in the
same way as to $U (K) $.  This gives the required homeomorphism of
the manifold $K$. \ep \medbreak

\section{On an extension of homeomorphisms from the torus to the
solid torus.}

The condition that the chosen curves are  homotopic for the
isomorphism of the schemes is not constructive, --- if the
homotopy condition is omitted in the definition of the isomorphism
of the schemes, then there are infinitely many different
isomorphisms of the schemes. Thus, it is impossible find among
them the one that satisfies the homotopy condition in a finite
number of steps. On the other hand, there  are infinitely many
different ways to construct, in relation to other curves on the
torus, a chosen cycle in the general case.  To overcome these
difficulties, we introduce a $ \tau $-invariant and replace the
homotopy condition  imposed on chosen cycles by the condition that
they have equal $ \tau $-invariants, which is easy to check.

Let $T $ be the boundary of a round handle and $G=T \cap S
^{(1)}$, where $S ^{(1)} $ is the $1$-skeleton of the stratified
set. We fix an orientation on $T $. Let $Z = \{ z _{i} \}$ be the
set of non-self-intersecting cycles of the graph $G$. Set an
orientation on each of these cycles and on a chosen cycle $w $ in
$T $. Define a correspondence between each oriented cycle $z _{i}
$ and an integer $ \alpha _{i} $ that equals the algebraic number
of points where it intersects $w $.  Set $ \beta _{i} $ to be
equal to the algebraic number of points in which $z _{i} $
intersects $v$, taken modulo $ \alpha _{i} $, where $v $ is a
parallel on $T$.  Thus, we have a pair of maps, $ \alpha , \beta :
Z \to \textbf{Z}$. If an isomorphism of the graphs is given, $
\phi : G _{{ 1}} \to G _{{ 2}} $, then it induces a pair of maps $
\phi^* ( \alpha _{1} ), \phi^* ( \beta _{1} ): Z _{2} \to \textbf
{Z}$ from the pair of maps $ \alpha _{{ 1}} , \beta _{{ 1}} : Z
_{{ 1}} \to \textbf{Z}$.

\begin{lemma}  Let $ T _{{ 1}} $ and $T _{{2}} $ be
boundaries of solid tori, $G_{i} \subset T _{i} , \ i = 1,2$, be
the embedded graphs. Suppose that the isomorphism of  the graphs,
$ \phi : G _ {{ 1}} \to G _{{ 2}} $, is   extended to a
homeomorphism  that preserves the orientation of the tori.   Then
it can be extended (probably by using another homeomorphism of the
tori) to a homeomorphism of the solid   tori if and only if $
\phi^* ( \alpha _{1})=\alpha_{2}, \phi^*(\beta_{1}) = \beta _{2}$.
\end{lemma}

\proc{Proof.} The necessity is obvious. To prove sufficiency, we
construct homeomorphisms of the tori which map chosen cycles into
the chosen cycles. There are three essentially different cases
that are possible:

1) all cycles are homotopic to $0$.  This is equivalent to the
graphs being homotopic to $0$.  Then, for every $i = 1,2$, one of
the components of $T _{i} \backslash G _{i} $ has genus $1$. Fix a
chosen cycle in each  component. Then it is easy to construct a
homeomorphism between these components which would map the chosen
cycles into the chosen cycles (cutting a component along a chosen
cycle we then need to extend the homeomorphism from the boundary
of the sphere with holes to its interior). On other components of
$T _{i} \backslash G _{i} $, we define the homeomorphisms that
satisfy the  conditions of the lemma.  Since all these
homeomorphisms coincide on the edges, they generate a
homeomorphism of the torus mapping the chosen cycle into the
chosen cycle.  The existence of such a homeomorphism is a
sufficient condition for existence of extensions of the
homeomorphism inside a solid torus.

2) There exist cycles that are not homotopic to $0$ but all of
them are homotopic to each other. This means that there also
exists a closed curve $ \gamma $ that is homotopic to these cycles
and not intersecting the graphs.  The pair $ \alpha , \beta $ is
then defined by the number of points where $ \gamma $ intersects
the chosen cycle and the parallel. Since the image of a chosen
cycle of $T _{{ 1}} $ and a chosen cycle on $T _{{ 2}} $ have the
same intersection numbers under the homeomorphism of the tori, $h
: T _ {{ 1}} \to T _{{ 2}} $, they differ from each other, up to a
homotopy, in a Dehn twist along $ \gamma $.  Since this Dehn twist
is a homeomorphism leaving the graphs fixed, we have a
homeomorphism of toruses mapping a chosen cycle into a chosen
cycle.  Hence, it can be extended to a homeomorphism of solid
toruses.

3) There are two or more cycles that are not homotopic to each
other and not homotopic to $0$. It is sufficient to consider two
such cycles and the corresponding numbers $ \alpha _{i} $. Using
this, the chosen cycle can be reconstructed uniquely up to a
homotopy.  Since the numbers $ \alpha _{i} $ of the corresponding
cycles are same,  the homeomorphism of the
  tori maps the chosen cycle into a chosen cycle,
which was to be shown. \ep \medbreak

\proc{Remark.} It is not necessary to calculate the map $ \beta $
in the third case. \medbreak

 There is a   principal difference
between the first two cases and the third one. It consists in
extending the homeomorphism of tori to a homeomorphism of the
solid in the third case as oppose to replacing the homeomorphism
of the tori with another one in first two cases  and then
extending it to solid toruses, under the conditions of the lemma.

In the first case, the graphs divides the torus into parts, one of
which has genus $1$ (a torus with holes). Denote it by $L $.  Then
$L $ is simultaneously a part of the boundary of a round
$0$-handle and a part of the boundary of a round $2$-handle (for
others types of handles, their boundaries consist of strata of
genus $0$). For a meridian of a round $2$-handle with respect to a
round $0$-handle, define the pair of numbers $( \alpha , \beta )$
as above. Then this pair is an invariant of the flow. For
convenience, we assume that $ \alpha =0$ and $ \beta =0$ for other
types of intersections of round $0$- and $2$-handles.

In the second case, it is possible to make a Dehn twist along $
\gamma $ before extending the homeomorphism of the torus. In such
a case, it is possible  that $ \gamma $ will not be homotopic to a
$0$-curve in the other torus, which means  that the Dehn twist
along it changes the homeomorphism of this torus, which is not
desirable. Therefore, it is necessary to introduce an additional
invariant that would indicate this situation. If $ \gamma $ is
homotopic to $0$ or to a chosen curve in the second torus, then
the Dehn twist along it does not influence the possibility of
extending the homeomorphism of the second torus. Hence, we will
consider the situations in which the algebraic numbers of
intersection points of $ \gamma $ with both meridians are not
equal to $0$. We choose, on each of the two tori, a closed curve $
\omega $ that intersects $ \gamma $ in one point and such that,
for the homotopy class of a meridian $m $, we have $[m] $ = $ k
[\gamma ] + l [\omega ]$, where $0 < k < \vert l\vert $. Then the
condition that the meridians are homotopic can be replaced by the
condition that the curves $ \omega $ are homotopic and the pairs
of numbers $(k , l )$ are equal. Note that $m =\omega $ for round
$1$-handles.

If there are two homotopic curves $ \gamma $ and $ \gamma ' $ in
different components of $T\backslash G $ on the boundary of the
round handle, then simultaneous Dehn twists along $ \gamma $ and
along $ \gamma ' $ with opposite orientations does not change the
homotopy class of $m $ and $ \omega $. Furthermore, the algebraic
numbers of intersection points of $ \omega $ with $ \omega $ in
other toruses are changed, but the sum of these numbers remains
constant.  Consider the graph $L $ whose vertexes correspond to
the usual  and round handles  containing domains $D$ that are
homeomorphic to rings  in case 2. The edges of the graph
correspond to these domains.  The ordering of handles defines a
natural orientation of the edges of $L $. To each edge of $L$, we
put into the correspondence a number $ \mu $ that equals the
algebraic number of intersection points of the curves $ \omega $
in the  corresponding domain, or $ \infty $, if the curve $ \gamma
$ is homotopic to $0$ in the boundary of one of the two handles to
which it belongs.  Here,  we take the curve $ \omega $ to be the
first if the order of its handle is less.  The operation  of
simultaneous twisting on the curves $ \gamma $ and $ \gamma ' $
generate an equivalence ration on the set of numbers $ \mu $,
which will be considered in more detail in the following section.

We fix one point by one component of  $\partial D$. We demand that
curves $\omega$ intersect $\partial D$ only in these points. For
two curves $\omega_1$ and $\omega_2$ that intersect $\partial D$
we use its isotopy to curves with such property : if we chance
parallel orientations on components of $\partial D$ and
orientation of $\omega_1 \bigcup \omega_2 \bigcap D$ then the
orientations are same in one points of $\omega_1 \bigcup \omega_2
\bigcap \partial D$ and opposite in another 



\section{Equivalence of framed graphs}

In this section, we consider oriented graphs in which one of the
ends of each edge is a source or a sink and we call them
MS-graphs. Except for the situation described above about the
graph $L$, such graphs are a   main invariant of a Morse--Smale
flow without closed orbits on closed surfaces. Their vertexes are
fixed points and their edges are separatrices. Vertexes of the
MS-graph, which are neither sources nor  sinks, are called
saddles.  In the graph $L $, each round $0$-handle corresponds to
a source of $L$ and round $2$-handles correspond to a sink, and
each saddle of $L$ corresponds to a round $1$-handle (a round
$1$-handle can also correspond to a source or a sink of $L$).

An oriented graph are called framed, if each edge is put into a
correspondence with an integer number or $ \infty $. This number
or $\infty$ are called a framing of the edge and the set of them
is called a framing of the graph. Thus, a framing of a graph is a
map of the set of edges into $\textbf{Z} \cup \ \{\infty\}$. A
framed graph is similar to a graph with colored edges.

Two framings of an MS-graph are called equivalent, if one of them
can be obtained from another by a sequence of the following
operations:

1) Simultaneous addition of an integer $k $ to the framing of two
edges of the graph, provided they have the property that the
origin of one is an endpoint of the other  and, hence, this point
is a saddle. Edges with such property will be called sequential.

2) Addition of an integer $k $ to the framing of one edge and $- k
$ to the framing of the incident edge, provided that these edges
are not sequential.

\begin{lemma}   If each of the two framings of  an MS-graph
contains edges with the framing $ \infty $, then they are
equivalent if and only if their sets of edges with  framing $
\infty $ coincide.
\end{lemma}

\proc{Proof.} Necessity. The operations 1) and 2) do not change
the set of edges with framing $ \infty $.

Sufficiency. Let us show that, under conditions of the lemma, any
edge with a finite framing $n $ can be changed to any integer $k $
without changing the framings of the remaining edges.  Indeed, in
view of the connectedness of the graph, there is a path $\{e _ {{
1}} , \ldots, e _{m} \}$ that starts at this edge and ends at an
edge with framing $ \infty $. Applying subsequently the operations
1) or 2) with numbers $k $ and $- k$ to pairs of the edges $(e _{{
1}} , e _{{ 2}} )$, $( e _ {{ 2}} , e _{{ 3}} )$, \ldots, $(e _{m
-{ 1}} , e _ {m} )$ we obtain the necessity statement. \ep
\medbreak

The main aim of this section is to establish criteria for  an
equivalence of framed MS-graphs. Assume that the graphs are
connected and do not contain edges with framing $ \infty $. We
subdivide all MS-graphs into 3 types and formulate criteria for
equivalence for each type.

1. The graphs of the first type are MS-graphs without saddles.
They allow only the second operation that, as easily  seen, does
not change the sum of the framings of all the edges.

\begin{lemma}  Two framings of the graph of  the first
type are equivalent if and only if they have the same sum  of the
framings of all edges.
\end{lemma}

\proc{Proof.} The necessity is obvious.

Sufficiency. Fix an edge $e $. Similarly to Lemma 1, for any edge
$e _{{ 1}} $, by choosing a path between this edge and the edge $e
$, we can  change the framing of $e _{{ 1}} $ by an arbitrary
integer $n $ so that the framing of $e $ will be changed by $-n$
with framings of the other edges preserved.  Then using operation
2) for the first framing graph we successively change the framing
of each edge to the necessary value modifying the value of the
frame of the edge $e $ by the used difference. As a result, we
have that it is possible to obtain, from the first framing, a
framing that coincides with second one except for, probably, the
edge $e $.  Since the sum of  the framings of the edges are the
same, these framings also coincide on $e $. \ep \medbreak

\begin{corr}   A framing of a graph of the first type
is equivalent  to a framing of a graph in which all edges have
framing $0$   save for a fixed edge $e$.
\end{corr}

2.  MS-graphs $G$ belong to the second type if there is a cycle
(probably nonoriented) that has an odd number of saddles ($G$ is
considered as an oriented graph).  It is clear that the saddles of
the cycle are saddles of the MS-graph (the inverse, generally, is
not true).  Existence of such a cycle is equivalent to existence
of a self non-intersecting cycle that has an odd number of
saddles.

\begin{lemma}   Two framed graphs of the
  second type are equivalent if and only if the sums of
  their framings of all edges are equal modulo 2.
\end{lemma}

\proc{Proof.} \textit{ Necessity.} Each of the operations 1) and
2) does not change the sum of framings of all edges modulo 2.

\textit{Sufficiency.} Let $e $ be an edge of a self
non-intersected cycle $ \omega = \{e , e _{ 1} , \ldots, e _{m}
\}$ with an odd number of saddles. As in the previous lemmas, by
using operations 1) and 2), we can pass from the first framing to
a framing that coincides with the second one on all edges except
for the edge $e $, for which the difference of these framings is
even.  A
  subsequent modification of the framings by $1$ or $-1$
  applying admissible operations to the
pairs of edges $\{e , e _{ 1} \}, \{e _{{ 1}} , e _{{ 2}} \},
\ldots, \{e _{m}, e \}$ changes the framing of the edge $e $ by
$2$ leaving the framings of other edges unchanged. Thus, we can
make the framing of $e $ equal to the second framing of the graph.
\ep \medbreak

3. The remaining MS-graphs are of the third type. We cut the
graphs into saddles. Since there are no cycles of the second type,
the components of the obtained graph can be subdivided into two
groups in such
  a way that, after the slitting, any two
subsequent edges belong to components from distinct groups.  One
of groups will be called the first one, another one the second.
The sum of framings of the edges in a group will be called the
total framing of the group.

\begin{lemma}  Two framings of the third
  type are equivalent if and only if they have the same
  differences of total framings of the first and second groups.
\end{lemma}

\proc{Proof.} \textit{Necessity.} If two adjacent edges belong to
one group, then the operation 2) can be applied to them. This
operation does not change the sum of the framings of the group. If
the adjacent edges belong to different groups, then the operation
1) can be applied to them; this operation does not change the
difference of the total framings of the groups.

\textit{Sufficiency.} Similarly to the previous lemmas, we can use
admissible operation to
  obtain, from the first framing, a framing
that coincides with the second one, except for, possibly, one
edge.  The conditions of the lemma guarantee that the framings of
this edge coincide. \ep \medbreak

In the case of disconnected graphs, the framings are equivalent if
the framings are equivalent for each connected component.

\section{The main theorem}

In the previous section, we have constructed a set of invariants
that indicate a possibility of extending homeomorphisms from the
boundaries of round handles to their interior. We call this set a
$ \tau $-invariant. Thus, a $ \tau $-invariant consists of the
  following data:

1) pairs of numbers $( \alpha $, $ \beta )$ for each pair of round
$0$- and $2$-handles that intersect on a surface of the genus $1$
with the orientation of this surface as in the boundary of the
round $0$-handle;

2) pairs of numbers $( \alpha $, $ \beta )$ for meridians of round
handles of type $2$, pairs $(k , l) $ for the curves $ \omega $,
and the equivalence class of the graph $L $;

3) the numbers $ \alpha _{i} $ for every non-self-intersecting
cycle on the boundaries of the round handles of type $3$.

\medskip

\begin{thm} Let $X$ and $X'$ be  Morse--Smale flows on a
$3$-manifold $M$. There is a topological equivalence of two
Morse--Smale flows $X$  and $X'$ preserving the order of critical
elements if and only if there is an isomorphism of their schemes
preserving the $ \tau $-invariant, chosen curves, and their
orientations.
\end{thm}

\proc{Note.} As opposed to Theorem 1, this theorem is
constructive, that is, it allows, in a finite time, to make a test
on topological equivalence of two flows.

\proc{Proof.}  \textit{Necessity}. It follows from the
construction and the results of the previous paragraph.

\textit{Sufficiency}. As in the proof of Theorem 1, we extend the
homeomorphism of the complexes $K $, possibly fixing it, to a
homeomorphism of the manifold $K^3$ without neighborhoods of
closed orbits of index $0$ and $2$.  Let us show that, if the $
\tau $-invariants are equivalent, then the homeomorphisms from the
tori that are boundaries of neighborhoods of closed orbits of
indices $0$ and $2$ can be extended inside.

In case 1 of the definition of the $ \tau $-invariant, by using
the map $g $, the homeomorphism between the tori can be changed as
in Theorem~1, since the images of the boundaries of the disks,
which  intersect the closed orbit transversally in a single point,
are homotopic to curves with the same properties. Here we use the
fact that the  corresponding numbers (the $ \tau $ invariants) are
equal. This means that these boundaries differ by Dehn twists
along the line of intersection of the stable manifolds of the
lower orbit and the torus.  The Dehn twist is performed by
modifying the homeomorphism of the diagram with a  contraction or
expansion that is isotopic  to the intersection of unstable
manifolds of this lower orbit and the surface $\Phi $.  Thus, as
in the cases 2 and 3, boundaries of disks are mapped into curves
that are homotopic to boundaries of disks, hence they  are also
boundaries such disks.

Using the first return map, extend the homeomorphisms of the tori
to disks and then, along the trajectories, to the  whole solid
torus.  This  gives a needed homeomorphism between the manifolds.
\ep \medbreak

\section{The distinguishing graph of the  flow }

Now we construct an invariant of the diagrams, which is complete
up to a (strict) isomorphism. It consists of a graph with
additional information.

The graphs $G $ is formed by the curves $u _{i} $ and $v _{j} $,
the curves in the intersections $\Phi \cap F \cap T _{i} $, and
the curves that define the $ \tau $-cycles. The vertexes of the
graph $G $ are intersections points of these curves; the edges are
arcs between them.  If a curve does not intersect other curves,
then it forms a loop on the graph (one vertex and one edge with
the ends in it).  We denote the vertexes by $A _{i} $ and the
edges by $B _{i} $. Fix an arbitrary orientation on the edges on
which it is not given.

By cutting the complex $K $ along the graph $G$, we obtain a set
of surfaces $F _{i} $ with boundaries on which we fix
orientations.  For each component of the boundary, passing it
along its orientation that agrees with the orientation of the
surface $F _{i} $, we write a word that consists of the letters $B
_{j} ^{{ \pm}{ 1}{ }} $ corresponding to the edges passed. The
letter has power $+1$ if the orientation of the corresponding arc
coincides with the orientation of the circle, and $-1$ in the
opposite case.  Two words are called equivalent if one can be
obtained from the other by a cyclic permutation of the letters.
This corresponds to another choice for the beginning of the walk
on the circle. The words are called inverse if one  is obtained
from the other by rewriting the letters in the inverse order
reversing the  powers and, possibly, making a cyclic permutation.
This corresponds to a walk on the circle in the direction opposite
to the orientation.

For each surface $F _{i} $, we form a list that consists of the
following:

1) the number $n _{i} $ that equals the genus of the surface $F
_{i} $, if the surface is oriented, and $ -n _{i} $, otherwise;

2) the words that are written when passing the boundary of the
surface along the orientation.

Such two lists are called equivalent if the numbers $n _{i} $ are
equal and there is a bijection between the words such that the
corresponding words are equivalent or inverse to
  each other.

Thus, we have constructed a set of lists of words (SLW) in such a
way that each list corresponds to one surface $F _{i} $.  Two such
sets are called equivalent, if there is a bijective correspondence
between the lists such that the corresponding lists are
equivalent.

Form a set of lists of words for each boundary of the round handle
and, also, a subset of the incoming and outgoing sets for the
round 1-handles.

We also make lists for the
  following curves:

1) the curves $u _{i} $ broken into pairs;

2) the curves $v _{j} $ broken into pairs;

3) the chosen cycles for the mean orbits;

4) the cycles and circles included in the $ \tau $-invariant in
case 2;

5) the curves included in the $ \tau $-invariant in case 3.

For the $ \tau $-invariant in case 1, we assign the corresponding
number to the regions that lie on the torus.

A graph $G$ with SLW and lists 1) -5) is called a distinguishing
graph of the flow.

Two distinguishing graphs are called equivalent if there is an
isomorphism of the graphs mapping the SLW of first graph into an
SLW that is equivalent to the SLW of the second graph and the
lists 1) -5) into equivalent lists preserving the
$\tau$-invariant.

\begin{thm}  Two Morse--Smale flows on  a
closed 3-manifold are topologically  equivalent if and only if
they have minimal diagrams  whose distinguishing graphs are
equivalent.
\end{thm}

\proc{Proof.} Taking into consideration the reasoning used in
Section 7, it is necessary to prove that the equivalence of the
distinguishing graphs is equivalent to stable equivalence of the
stably minimal diagrams. The proof of this is similar to the proof
of Theorem 1 in [8] with replacing the surface by the complex $K
$. \ep \medbreak

\section{Examples}

Consider Morse--Smale flows with one closed 1-orbit and the
minimal number of other critical elements on $S^3$.

First consider a flow with three closed orbits and no
singularities.  A closed 1-orbit can have a trivial or a twisted
neighborhood. In the case of a trivial neighborhood, the round
$1$-handle is attached to a round $1$-handle by two annuli such
that the core of one of them is homotopic to $0$ and another one
is not. 
So the distinguishing graph consists of 4
loops, $a, b, c, d$. SLW = \{$L_1$=\{0, $a, b^{-1}$\}, $L_2$=\{0,
$c$, $d^{-1}$\}, $L_3$=\{0, $a^{-1}$, $b$, $c^{-1}$ \}, $L_4$= \{
0, $d$ \}, $L_5$=\{0, $a$, $d^{-1}$\}, $L_6$=\{0, $b$, $c^{-1}$ \}
\}.  The  boundary of the round handles are $T_0=\{L_1, L_2, L_3,
L_4 \} $, $T_1=\{ \{L_1,L_2\}, \{L_5, L_6\}\} $, $T_2=\{L_3, L_4,
L_5, L_6\} $.  The $\tau$-invariant is non-trivial for the
intersection of $b$ with the meridians of the round $0$- and
$2$-handles.  The algebraic number of such intersections is equal
to $\pm 1$. Hence, $4$ cases are possible.



Thus, there are $4$ different flows. One can be obtained from
another be replacing the orientation of the closed $0$- or
$2$-orbit and the flow in its neighborhood.

In the case of a twisted $1$-orbit, the round $1$-handles are
regular neighborhoods of $2n+1$ $\pi$-twisted bands in the
complement of the round $0$-handle 
The
distinguishing graph is the loop $a, b$. SLW = $\{L_1= \{ 0, a,
b^{-1} \}$, $L_2= \{ 0, a, b^{-1} \}$, $L_3= \{ 0, a, b^{-1} \}
\}$.  The algebraic number of the intersection of the meridian of
the round $0$-handles and $a$ is equal to $2n+1$. It is a unique
non-trivial part of the $\tau$-invariant. Hence, there are
infinitely many such flows and they are given by $n$. The pairs of
numbers for $\omega$ of round 0-, 1- and 2-handles are equal to
$(1, 2n+1), (1,2)$ and $(0,1)$, correspondingly. Total framing of
$L$ can be equal 0 or 1 mod 2. It depend of the orientations of
the closed orbits.



If, in these examples, we replace each closed $0$- and $2$-orbit
by two singularities, then we obtain a flow with infinitely many
singular trajectories.

\proc{Acknowledgements.} I would like to thank Vladimir Sharko,
Yevgeniy Polulyah, Sergey Maksimenko and Igor Vlasenko for talks
promoting the appearance of this paper.

\end{document}